# The design of courier transportation networks with a nonlinear zero-one programming model

Boliang Lin*

School of Traffic and Transportation, Beijing Jiaotong University, Beijing 100044, People's Republic of China

**Abstract**: The courier industry is one of the most essential parts of the modern logistics. Meanwhile, the courier transportation network is one of the most important infrastructures of courier enterprises to take participate into operation. This paper presents a combinatorial optimization model for the courier transportation network design problem, and the aim of network optimization is to determine the transportation organization mode for each courier flow. A nonlinear zero-one integer programming model is formulated to describe the problem, the objective function of the model is intended to minimize the total cost including the accumulation cost, the transportation cost and the transfer cost. Also, we take transportation modes and types of transport carriers into account in the objective function. The constraints of the combinatorial optimization model contain capacities of transfer and sorting centers and delivery dates predefined.

## 1 Introduction

The rapid rise of e-commerce stimulated by the internet, the increase of courier demand caused by the network shopping and the development of various transportation modes have promoted the development of the courier industry. The courier industry has experienced fast development in China in the past decade. Take China as an example, in 2013, the growth rate of courier business reached the highest level in history which is 61.6%. According to the recent data, the courier business volume reached 31.28 billion pieces with a growth rate of 51.4%, which was 6.5 times the GDP growth rate of China in 2016. As the development growth of the economy tends to stabilize, the development speed of the courier industry also tends to be stable. In this situation, the optimization of the courier transportation network can promote the development of the courier industry further, and improve the economic development in turn.

Meanwhile, as the carriers of courier transportation, road transportation, railway transportation and air transportation have achieved great progress and proposed brand new courier transportation services to meet the courier transportation demand. Continuous competition and cooperation among those three transportation modes have promoted the improvement of China's transportation system and the development of courier industry.

---

* E-mail address: bllin@bjtu.edu.cn (B.-L. Lin).

Currently, sizeable Chinese courier enterprises include: S.F. Express, YTO Express, ZTO Express, STO Express, and Yunda Express, etc. The development scales of China's five largest courier enterprises by the end of 2016 are shown in Table 1.

Table 1. Five major courier enterprises' development scales as of the end of 2016[1]

|  | S.F.[2] | YTO[3] | ZTO[4] | STO[5] | Yunda[6] |
|---|---|---|---|---|---|
| Urban coverage above county-level | 91.6% | 96.1% | 96.6% | 96.4% | 95.0% |
| Terminal outlets | 13000 | 37713 | 26000 | 20000 | 20000 |
| Main transportation lines on Land | 9600 | 3475 | 1980 | - | 4200 |
| Transshipment centers | 272 | 62 | 69 | 48 | 55 |
| Self-operated vehicles | 15000 | 36000 | 2930 | - | - |
| Self-operated freighters | 36 | 5 | - | - | - |

Optimization of the network has a deep influence on the timeliness and costs. This paper established a nonlinear zero-one integer programming model to optimize the courier delivery network so that courier can be delivered to customers at a lower cost within the specified delivery time.

The remainder of the paper is organized as follows. Section 2 reviews relevant literature related to the courier transportation network and its optimization. Section 3 introduces the problem we try to solve in this paper. Section 4 presents the notations definition and modeling assumptions of the optimization model, and also gives a mathematical programming model for the problem. Section 5 concludes the paper.

## 2 Literature review

Before studying the design problem of the courier transportation network, we should have a certain degree understanding of the structure and characteristics of the courier transportation network. For the study of the courier operational network, Aykin (1995) introduced a framework for the design of hub-and-spoke distribution network. Jeong et al. (2007) described the hub-and-spoke rail networks. Tan (2015) studied the topology structure of the courier network, the traffic spatial-temporal dynamics and the package delay distributions based on the logistics data.

The courier transportation network design and the optimization problem are our focuses in this paper. Barnhart et al. (2002) focused on a particular service network design application to determine cost minimizing routes and schedules with time windows. Schwind and kunkel (2010) presented the research aiming at a dynamic optimization of tour planning process in courier, express and parcel delivery networks. Alibeyg et al. (2016) presented a class of hub network design problems with profit-oriented objectives integrating several locational and network design decisions. Di et al. (2018) studied a new discrete network design problem for metropolitan areas.

In addition, multimodal transport is generally applied to courier operational network. The intermodal hub network design problem began with the pioneering work of Slack (1990), who analyzed the changes that were transforming intermodal transportation in the United Stated and Canada.

There're always time constraints in the network design problem. Iyer and Ratliff (1990) organized a guaranteed time distribution system. Kara (2011) focused on the minimization of the arrival time of

---

the last arrived item in cargo delivery systems.

The contribution of this paper is the consideration of accumulation process of courier flows, the limited capacity of facilities and the constraint of delivery time based on real life conditions. The accumulation time calculated in this paper is the accumulation parameter multiplied by the standard courier unit. The delivery time of courier flows consists of transportation time, transfer operation time and accumulation time.

## 3 Problem description

Here, we use Fig. 1 as an example to illustrate the courier delivery process. Courier at each pick-up station is sent to the terminal delivery points, then delivered to local distribution centers, courier at local distribution centers can be sent to the airports, railway stations or transfer and sorting centers, or can be sent directly to a certain local distribution center of terminal city. After courier arrives at the local distribution center in the city, it'll be delivered to the nearby terminal distribution nodes by small van-body trucks, then, sent to the community outlets by small transport carriers.

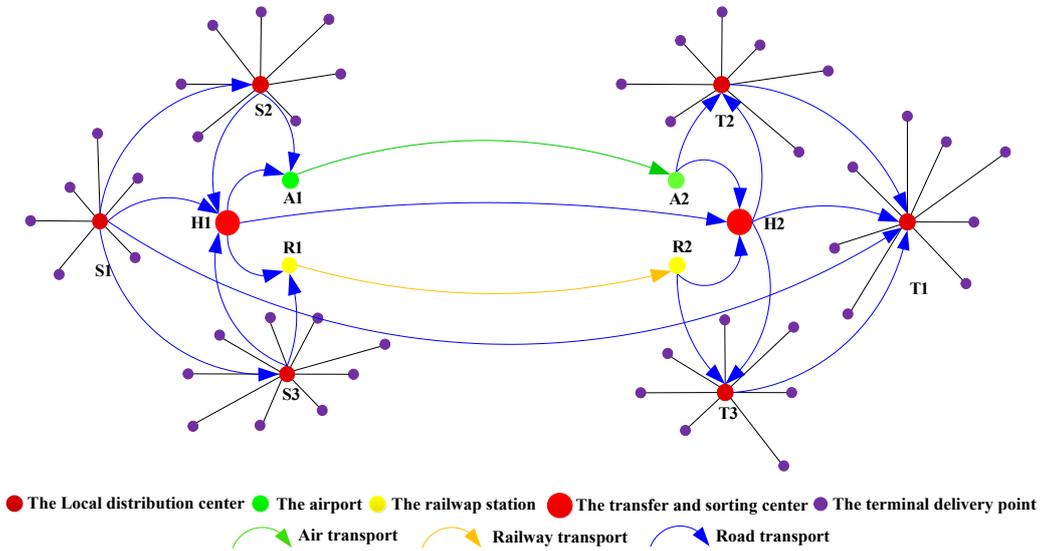

● The Local distribution center  ● The airport  ● The railwap station  ● The transfer and sorting center  ● The terminal delivery point
　Air transport　　Railway transport　　Road transport

**Fig. 1.** The courier transportation network

Here, we take the courier flow $N_{S1,T1}$ from $S1$ to $T1$ as an example to illustrate the design idea of the courier transportation network optimization. There are several possible transportation options from $S1$ to $T1$:

(1) $S1 \to T1$
(2) $S1 \to H1 \to T1$
(3) $S1 \to H2 \to T1$
(4) $S1 \to H1 \to H2 \to T1$
(5) $S1 \to A1 \to A2 \to T1$
(6) $S1 \to H1 \to A1 \to A2 \to T1$
(7) $S1 \to A1 \to A2 \to H2 \to T1$
(8) $S1 \to H1 \to A1 \to A2 \to H2 \to T1$
(9) $S1 \to R1 \to R2 \to T1$
(10) $S1 \to H1 \to R1 \to R2 \to T1$
(11) $S1 \to R1 \to R2 \to H2 \to T1$
(12) $S1 \to H1 \to R1 \to R2 \to H2 \to T1$

We need to consider the matter of which point pairs provide direct transportation services and what type of transfer option is used for each courier flow. The aim of optimization is to minimize the total cost of the whole courier transportation network with real-world constraints like delivery time, capability of transfer nodes etc., which constitutes the combinatorial optimization problem of courier transportation network design solved by this paper.

## 4. Mathematical formulation

This section will derive a mathematical model for the courier transportation network optimization problem according the above description. In order to facilitate readers' understanding of the model, here we first list all the notations involved in the modeling process (including sets, parameters and decision variables) as follows:

### 4.1 notations

**Sets**

| | | |
|---|---|---|
| $\mathbb{S}^{node}$ | | The set of all nodes in the courier transportation network. |
| $\mathbb{S}^{shpt}$ | | The set of original courier flows $N_{ij}$. |
| $\mathbb{S}_{st}^{chain}$ | | The transfer chains of courier flows which consist of origins, destinations and transfer nodes passed across by the courier flow, $\mathbb{S}_{st}^{chain} = \{k_1, k_2, \cdots, k_n\}$. |
| $\mathbb{P}(i,j)$ | | The set of potential transfer nodes (not including node $i$ and node $j$, the node $j$ in this paper refers specifically to local distribution centers of the terminal city and node $i$ can be any type of nodes) that may be passed through by the courier flow $f_{ij}$. |

**Parameters**

| | |
|---|---|
| $\lambda$ | The time value used to unify the dimension of each factor in the objective function. |
| $\mathbb{V}_{ij}^{size}$ | The volume of the standard courier unit loaded by a transport carrier (such as a truck of a certain tonnage capacity) if the direct transportation is available between node $i$ and node $j$. |
| $\mathbb{C}_i^{acum}$ | The accumulation parameter of node $i$ covering all the factors affecting the accumulation time except for the transport unit. |
| $\mathbb{C}_{ij}^{trans}$ | The transportation costs of a single transport unit from node $i$ to node $j$. |
| $R_k^{sort}$ | The daily sorting quantity of the sorting center $k$. |
| $\mathbb{C}_k^{opr}$ | The operation costs of sorting operation unit volume of sorting center $k$. |
| $N_{ij}$ | The daily original courier flows (in standard courier units) which origin at local distribution center $i$ and are destined to terminal local distribution center $j$ |
| $f_{ij}$ | The daily average courier flows from node $i$ to node $j$, including the original demand $N_{ij}$ at node $i$ and other flows sorted and transferred at node $i$ arriving at local distribution center $j$. |
| $F_{ij}^{serv}$ | The service flows from node $i$ to node $j$. |
| $\phi_{ij}^{serv}$ | The daily needed frequency of transport carriers on the direct service arc $i \to j$. |
| $\mathbb{H}_k^{xfer}$ | The transfer capacity at node $k$. |
| $t_{st}^{trans}$ | The travel time of the direct transportation from node $s$ to node $t$. |
| $t_k^{opr}$ | The transfer operation time at node $k$. |
| $t_{kv}^{\text{freq}}$ | The waiting time at node $k$ of the accumulated transferred courier in the direction from node $k$ to node $v$. |
| $T_{st}$ | The delivery time limits from node $s$ to node $t$. |
| $\mathbb{P}_{ij}^{\text{next}}$ | The next transfer node after node $i$ on the transfer chain $S_{ij}^{\text{chain}}$. |
| $t_{ij}^{\text{hdway}}$ | The headway on the direct transportation arc $i \to j$, which is the departure intervals of transport carriers. |
| **Decision Variables** | |
| $x_{ij}^k$ | $x_{ij}^k = 1$ if the first transfer node is $k$ when the courier flow is transported from node $i$ to node $j$, and is equal to zero otherwise. |
| $y_{ij}$ | $y_{ij} = 1$ if the courier flow can be transported from node $i$ to node $j$ directly, and is zero otherwise. |

## 4.2 Basic assumptions of the mathematical model

**Assumption 1:** Properties of the courier are consistent (Usually, the courier can be divided into two types according to the delivery dates offered: the ordinary courier and the express, which can be split into two types of courier networks when processing, thus, we only consider one of them here.).

**Assumption 2:** The relatively low weight of the courier results in the volume, not the weight, which transport carriers first restrict. Because of the different sizes of courier, we use a cubic meter (or a pallet allowable volume) as a standard courier unit (which probably is composed of dozens or hundreds of the original courier).

**Assumption 3:** A courier flow is not allowed to be split among multiple paths or transportation modes.

**Assumption 4:** The size of a batch of transportation volume $\mathbb{V}_{ij}^{\text{size}}$ (for example, for a 50-ton truck, standard courier unit could be $\mathbb{V}_{ij}^{\text{size}}=60$) is known if direct transportation service is provided between node $i$ and node $j$.

## 4.3 Description of constraints and intermediate decision variables

According to the above variable definition, the binary variable $x_{ij}^k=1$ represents the courier flow $f_{ij}$ is transported to node $k$ through certain transportation services from node $i$ (which could be the local distribution center, the airport, the railway station or the sorting center). The binary variable $y_{ij}=1$ expresses that the courier flow $f_{ij}$ is transported to the terminal distribution center $j$ directly. According to assumption 3, it's clear that there is a uniqueness condition for the direct and transfer transportation:

$$y_{ij} + \sum_{k \in P(i,j)} x_{ij}^k = 1 \tag{4}$$

Notice that it means the direct transportation service is available between node $i$ to node $j$ when $x_{ij}^k=1$ (where direct refers to no changes in the transportation mode and no sorting operations on the way, and it'll still be the direct transportation between node $i$ to node $j$ from a courier enterprise's perspective if the courier enterprise consigns shipment to the airport, and then the air transport department delivers shipment to node $k$ through several transfers in transit). Thus, there is a variable correlation constraint:

$$x_{ij}^k \leq y_{ik} \quad \forall i,j \in S^{\text{node}}, k \in P(i,j) \tag{5}$$

The variable $f_{ij}$ is an intermediate variable according to definition, whose specific expression is a recursive formula with the following form:

$$f_{ij} = N_{ij} + \sum_{s \in S^{node}} f_{sj} x_{sj}^{i} \quad \forall i, j \in S^{node} \tag{6}$$

If the direct transportation service is available between node $i$ and node $j$, and the daily direct volume is $F_{ij}^{serv}$, it is clear that:

$$F_{ij}^{serv}(X) = f_{ij} + \sum_{t \in S^{node}} f_{it} x_{it}^{j} \quad \forall i, j \in S^{node} \tag{7}$$

It's easy to find out that the direct volume here contains the transferred courier from nodes before node $i$ and courier flows whose destinations are node $j$ or other nodes behind node $j$. When the capacity of a transport carrier which serves on a direct transportation arc is $\mathbb{V}_{ij}^{size}$, the daily needed frequency of transport carriers on this direct transportation arc will be:

$$\phi_{ij}^{serv} = F_{ij}^{serv} / \mathbb{V}_{ij}^{size} \tag{8}$$

The courier flow which arrives at the distribution center and the sorting center is random, which makes the average waiting time of a courier unit that is on the direct transportation arc from node $i$ to node $j$ is the half of the headway.

$$t_{ij}^{freq} = \frac{1}{2} t_{ij}^{hdway} = \frac{24}{\phi_{ij}^{serv}} = \frac{12}{\phi_{ij}^{serv}} = \frac{12 \mathbb{V}_{ij}^{size}}{F_{ij}^{serv}} \tag{9}$$

For example, we analyze the waiting time (excluding the operation time) for the accumulation of courier flows at a sorting center. Assume that there are four trucks operated on the direct transportation arc $i \to j$, which is $\phi_{ij}^{serv} = 4$, making the headway is 12-hour, that is, the average frequency delay (the waiting time for the accumulation) of the courier flow in this circumstance is $t_{ij}^{freq} = \frac{1}{2} t_{ij}^{hdway} = 6$ hours.

Notice that the daily volume on the direct transportation arc $i \to j$ is $F_{ij}^{serv}$, thus, if the direct transportation service is offered on the arc $i \to j$, the total frequency delay of courier flows which is on this arc would be:

$$F_{ij}^{serv} t_{ij}^{freq} = 12 \mathbb{V}_{ij}^{size} \tag{10}$$

Considering that courier flows are not always arriving evenly, there could exists accumulation interruptions.

Therefore, the total frequency delay of courier flows should be expressed more strictly as follows:

$$F_{ij}^{serv} t_{ij}^{freq} = \mathbb{C}_{i}^{acum} \mathbb{V}_{ij}^{size} \tag{11}$$

Here, $\mathbb{C}_{i}^{acum}$ is less than 12 in general. In a real transportation process, courier flows don't arrive in a balanced manner usually, and its value is also affected by the size of a station and other factors, the specific value may generally be recommended to be 10 to 11.5.

If node $k$ is a sorting center, according to the definition of decision $x_{ij}^{k}$ and intermediate variable $f_{ij}$, we can conclude the daily sorting quantity of this sorting center would be:

$$R_{k}^{sort}(X) = \sum_{i \in S^{node}} \sum_{j \in S^{node}} f_{ij} x_{ij}^{k} \tag{12}$$

The total time consumption of all parts in transit of the courier flow should be less than the delivery date $T_{st}$ of the courier flow $N_{st}$ if $T_{st}$ is offered. The total time consumption is discussed in the following two parts:

(1) The total transportation time of the direct transportation

The total transportation time of the direct transportation would be $t_{st}^{trans}$ when $y_{st}=1$, that is, courier flow $N_{st}$ is transported directly from original local distribution center $s$ to terminal local distribution center $t$.

(2) The total transportation time of the transfer transportation

The transportation time required is composed of the travel time on each transportation arc and the transfer time at each transfer node if the courier flow $N_{st}$ is transferred several times before reaching the terminal local distribution center $t$. The transfer chain for courier flow $N_{st}$ denoted as $\boldsymbol{S}_{st}^{chain}=\{k_1,k_2,\cdots,k_n\}$, where $k_1$ is the origin $s$ of the courier flow and $k_n$ is the destination $t$, $k_2$ would be the first transfer node which the courier flow $N_{st}$ passes through. Recording the next transfer node of courier flow $N_{st}$ as $\mathbb{P}_{ij}^{next}$, $\boldsymbol{S}_{st}^{chain}$ can be defined as follows:

$$\boldsymbol{S}_{st}^{chain}=\left\{k_1=s, k_2=\mathbb{P}_{k_1,t}^{next},\cdots,k_n=t=\mathbb{P}_{k_{n-1},t}^{next}\right\}$$

The recursive formula of transfer nodes can be expressed as:

$$\mathbb{P}_{ij}^{next} = jy_{ij} + \sum_{k\in\boldsymbol{P}(i,j)} kx_{ij}^k \tag{13}$$

According to the expression of the transfer chain, the total travel time of arcs should be:

$$\sum_{\overrightarrow{uv}\in\boldsymbol{S}_{st}^{chain},\overrightarrow{uv}\neq\overrightarrow{st}} t_{uv}^{trans}$$

In the above formula, $\overrightarrow{uv}\in\boldsymbol{S}_{st}^{chain}$ means $u\to v$ is one of the transportation arcs of the transfer chain $\boldsymbol{S}_{st}^{chain}$, the time delay on the nodes of the transfer chain $\boldsymbol{S}_{st}^{chain}$ is the sum of transfer time at each node. For the sorting center, the transfer time consists of operation time and accumulation time:

$$\sum_{\overrightarrow{kv}\in\boldsymbol{S}_{st}^{chain},k\neq s,t} \left(t_k^{opr}+t_{kv}^{freq}\right)$$

## 4.4 The nonlinear zero-one integer programming model

According to the above discussion, the optimization model of the courier transportation network design can be presented as:

$$\min Z = \sum_{i\in\boldsymbol{S}^{node}}\sum_{j\in\boldsymbol{S}^{node}}\left[\lambda\mathbb{C}_i^{acum}\mathbb{V}_{ij}^{size}+\phi_{ij}^{serv}(X)\mathbb{C}_{ij}^{trans}\right]y_{ij} + \sum_{k\in\boldsymbol{S}^{node}}R_k^{sort}(X)\left(\lambda t_k^{opr}+\mathbb{C}_k^{opr}\right) \tag{14}$$

s.t.

$$y_{ij}+\sum_{k\in\boldsymbol{P}(i,j)} x_{ij}^k =1 \qquad \forall i,j\in\boldsymbol{S}^{node} \tag{15}$$

$$x_{ij}^k \leq y_{ik} \qquad \forall i,j\in\boldsymbol{S}^{node}, k\in\boldsymbol{P}(i,j) \tag{16}$$

$$R_k^{sort}(X)\leq \mathbb{H}_k^{xfer} \qquad \forall k\in\boldsymbol{S}^{node} \tag{17}$$

$$y_{st}t_{st}^{trans}+\sum_{\overrightarrow{uv}\in\boldsymbol{S}_{st}^{chain},\overrightarrow{uv}\neq\overrightarrow{st}} t_{uv}^{trans}+\sum_{\overrightarrow{kv}\in\boldsymbol{S}_{st}^{chain},k\neq s,t}\left(t_k^{opr}+t_{kv}^{freq}\right)\leq T_{st} \qquad \forall N_{st}\in\boldsymbol{S}^{shpt} \tag{18}$$

$$y_{ij},x_{ij}^k\in\{1,0\} \qquad \forall i,j\in\boldsymbol{S}^{node}, k\in\boldsymbol{P}(i,j) \tag{19}$$

In the objective function (14), the first term accounts for the accumulation costs and transportation costs in the direct transportation, the second term calculates the transfer operation costs and transfer time costs of transfer transportation. Constraint (15) ensures that each courier can reach its destination and limits the transportation organization mode for courier that can only be selected for the direct transportation or the transfer transportation. Constraint (16) limits the connection between nodes:

$y_{ik}=0, x_{ij}^{k}=0$ means that the direct transportation is unavailable between node $i$ and node $k$, and $y_{ik}=1, x_{ij}^{k}=0$ or 1 means that node $k$ could be the first transfer node on the path from node $i$ to node $j$ when direct transportation is available between node $i$ and node $k$, that is, the courier flow from node $i$ to node $j$ can choose the direct transportation or the transfer transportation. Constraint (17) is the capacity constraint at the transfer node. Constraint (18) is the delivery time limitation, the delivery time includes transportation time, transfer operation time and accumulation time. Constraint (19) is the binary constraints on the decision variables.

## 5. Conclusion

In this paper, we developed a nonlinear zero-one programming model to optimize the courier network. We considered three transportation modes: road transport, air transport and railway transport, and different types of transport carriers. In this model, total courier delivery time includes transportation time, transfer time and accumulation time. This model determines that direct transportation service is offered between which node pairs and transfer transportation service should be provided at which nodes, aiming to minimize the total transportation cost containing the direct transportation cost of which are composed of the accumulation cost and the transportation cost and the transfer transportation cost of which consist the transfer operation cost and the transfer operation time cost.

For future research, we can discuss how to choose the type of transport carriers served on a courier transportation service arc since types of transport carriers are already known in this paper. Also, the number of transport carriers on the road should be restricted to avoid a congestion, which could be considered in the model.